%% file: main.tex
\documentclass[graybox]{svmult}

\usepackage{mathptmx}       
\usepackage{helvet}         
\usepackage{courier}        
\usepackage{type1cm}        
\usepackage{makeidx}         
\usepackage{graphicx}        
\usepackage{multicol}        
\usepackage[bottom]{footmisc}

\usepackage{amssymb}
\usepackage{amsmath}

\makeindex             

\begin{document}

\title*{Preconditioned space--time boundary element methods for the 
one-dimensional heat equation}
\titlerunning{Preconditioned space--time BEM for the one-dimensional 
heat equation}
\author{Stefan Dohr and Olaf Steinbach}
\institute{Stefan Dohr, Olaf Steinbach \at Institut f\"ur Angewandte
  Mathematik, TU Graz, Steyrergasse 30, 8010 Graz \at
  \email{stefan.dohr@tugraz.at, o.steinbach@tugraz.at}}
%
%
\maketitle

\newcommand{\norm}[1]{\left\lVert#1\right\rVert}
\newcommand{\abs}[1]{\ensuremath{\left\vert#1\right\vert}}
\newcommand{\traceH}{H^{1/4}(\Sigma)}
\newcommand{\dualtraceH}{H^{-1/4}(\Sigma)}

\section{Introduction}
Space--time discretization methods, see, e.g., \cite{steinbach}, 
became very popular in recent years,
due to their ability to drive adaptivity in space and time simultaneously,
and to use parallel iterative solution strategies for time--dependent problems.
But the solution of the global linear system requires the use of some
efficient preconditioner. 

In this note we describe a space--time boundary 
element discretization of the heat equation and an efficient and robust preconditioning strategy which is based on the use of boundary integral operators of opposite orders, but which requires a suitable stability condition for the boundary element spaces used for the discretization. We demonstrate the method for the simple spatially one-dimensional case. However, the presented results, particularly the stability analysis of the boundary element spaces, can be used to extend the method to the two- and three-dimensional problem \cite{dohr_niino_steinbach}. 

Let $\Omega = (a,b) \subset \mathbb{R}$, 
$\Gamma := \partial \Omega = \left\{a, b \right\}$ and $T > 0$. 
As a model problem we consider the Dirichlet boundary value problem for 
the heat equation,
\begin{equation}
\label{eq_model_problem}
\alpha \partial_{t} u - \upDelta_{x}u = 0 \; \textrm{in} \; 
Q:= \Omega \times (0,T), \;
u = g \; \textrm{on} \; \Sigma := \Gamma \times (0,T),\;
u = u_{0} \; \textrm{in} \; \Omega
\end{equation}
with the heat capacity constant $\alpha > 0$, the given initial datum 
$u_{0}$, and the boundary datum $g$. The solution of \eqref{eq_model_problem} 
can be expressed by using the representation formula for the 
heat equation \cite{costabel},
i.e. for $(x,t) \in Q$ we have
\begin{equation}
\label{eq_representation_formula}
\begin{aligned}
u(x,t) = &\int_{\Omega} U^{\star}(x-y, t) u_{0}(y) \D y + \frac{1}{\alpha} \int_{\Sigma} U^{\star}(x-y, t-s) \frac{\partial}{\partial n_{y}} u(y,s) \D s_{y} \D s\\
&- \frac{1}{\alpha} \int_{\Sigma} \frac{\partial}{\partial n_{y}} U^{\star}(x-y, t-s) g(y,s) \D s_{y} \D s,
\end{aligned}
\end{equation}
where $U^{\star}$ denotes the fundamental solution of the heat equation 
given by
\begin{equation*}
\label{eq_fundamental_solution}
U^{\star}(x-y,t-s) =
  \begin{cases}
   \displaystyle \left(\frac{\alpha}{4 \pi (t-s)}\right)^{1/2} \exp{\left(\frac{-\alpha|x-y|^{2}}{4 (t-s)}\right)}\;,  \quad &s < t, \\
   \displaystyle \; 0\;,        &\mathrm{else}.
  \end{cases}
\end{equation*}
Hence it suffices to determine the yet unknown Cauchy datum 
$\partial_{n}u_{|\Sigma}$ to compute the solution of \eqref{eq_model_problem}. 
It is well known \cite{lions_magenes_2} that 
for $u_{0} \in L^{2}(\Omega)$ and $g \in H^{1/2, 1/4}(\Sigma)$ the 
problem \eqref{eq_model_problem} has a unique solution
$u \in H^{1, 1/2}(Q, \alpha \partial_{t} - \upDelta_{x})$ 
with the anisotropic Sobolev space
\begin{equation*}
H^{1, 1/2}(Q, \alpha \partial_{t} - \upDelta_{x}) := \left\{u \in H^{1, 1/2}(Q): (\alpha \partial_{t} - \upDelta_{x})u \in L^{2}(Q) \right\}.
\end{equation*}
In the one-dimensional case the spatial component of the space--time 
boundary $\Sigma$ collapses to the points $\left\{a, b \right\}$ and 
therefore we can identify the anisotropic Sobolev spaces 
$H^{r, s}(\Sigma)$ with $H^{s}(\Sigma)$. The unknown density 
$w := \partial_{n}u_{|\Sigma} \in \dualtraceH$ can be found by applying 
the interior Dirichlet trace operator 
$\gamma_{0}^{\textrm{int}}: H^{1, 1/2}(Q) \rightarrow \traceH$ to the
representation formula \eqref{eq_representation_formula}, 
\begin{equation*}
g(x,t) = (M_{0}u_{0})(x,t) + (Vw)(x,t) + ((\frac{1}{2}I-K)g)(x,t) 
\quad \textrm{for } (x,t) \in \Sigma.
\end{equation*}
The initial potential $M_{0}: L^{2}(\Omega) \rightarrow H^{1/4}(\Sigma)$, 
the single layer boundary integral operator 
$V: H^{-1/4}(\Sigma) \rightarrow H^{1/4}(\Sigma)$, 
and the double layer boundary integral 
operator $\frac{1}{2}I-K: H^{1/4}(\Sigma) \rightarrow H^{1/4}(\Sigma)$ 
are obtained by composition of the potentials in 
\eqref{eq_representation_formula} with the Dirichlet trace operator 
$\gamma_{0}^{\textrm{int}}$, see, e.g., \cite{costabel, noon}.
In fact, we have to solve the variational formulation to find 
$w \in \dualtraceH$ such that
\begin{equation}
\label{eq_var_form_first_bie}
\langle Vw, \tau \rangle_{\Sigma} = 
\langle (\frac{1}{2}I + K)g, \tau \rangle_{\Sigma} - 
\langle M_{0}u_{0}, \tau \rangle_{\Sigma} \quad 
\textrm{for all } \tau \in \dualtraceH,
\end{equation}
where $\langle \cdot, \cdot \rangle_{\Sigma}$ denotes the duality pairing 
on $H^{1/4}(\Sigma) \times H^{-1/4}(\Sigma)$.
The single layer boundary integral operator $V$ is bounded and 
elliptic, i.e. there exists a constant $c_{1}^{V} > 0$ such that
\begin{equation*}
\langle Vw, w \rangle_{\Sigma} \geq c_{1}^{V} \norm{w}_{\dualtraceH}^{2} \quad 
\textrm{for all } w \in \dualtraceH.
\end{equation*}
Thus, the variational formulation 
\eqref{eq_var_form_first_bie} is uniquely solvable. When applying the 
Neumann trace operator 
$\gamma_{1}^{\textrm{int}} : H^{1, 1/2}(Q, \alpha \partial_{t} - \Delta_{x}) 
\to \dualtraceH$ to the representation formula
\eqref{eq_representation_formula}
we obtain the second boundary integral equation
\begin{equation*}
w(x,t) = (M_{1}u_{0})(x,t) + ((\frac{1}{2}I + K')w)(x,t) + (Dg)(x,t) \quad \textrm{for } (x,t) \in \Sigma
\end{equation*}
with the hypersingular boundary integral operator
$D : \traceH \to \dualtraceH$,
and with the adjoint double layer boundary integral operator
$K' : \dualtraceH \to \dualtraceH$. Moreover, $M_{1}: L^{2}(\Omega) \rightarrow \dualtraceH$. 

\section{Boundary element methods}
\begin{figure}[t]
\sidecaption[t]
\includegraphics[]{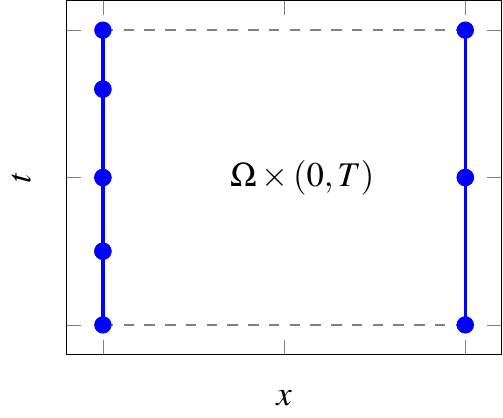}
\caption{Sample BE mesh. We consider an arbitrary decomposition of the 
space--time boundary $\Sigma$. Note that there is no time-stepping 
scheme involved}
\label{fig_mesh_1d}
\end{figure}

For the Galerkin boundary element discretization of the variational formulation
\eqref{eq_var_form_first_bie}
we consider a family $\left\{\Sigma_{N}\right\}_{N \in \mathbb{N}}$ of 
arbitrary decompositions of the space--time boundary $\Sigma$ into 
boundary elements $\sigma_{l}$, i.e. we have
\begin{equation*}
\overline{\Sigma}_{N} = \bigcup_{\ell=1}^{N} \overline{\sigma}_{\ell} \;.
\end{equation*}
In the one-dimensional case the boundary elements $\sigma_{\ell}$ are line 
segments in temporal direction with fixed spatial coordinate 
$x_{\ell} \in \left\{a, b\right\}$ as shown in Fig. \ref{fig_mesh_1d}. 
Let $(x_{\ell}, t_{\ell_{1}})$ and $(x_{\ell}, t_{\ell_{2}})$ be the nodes 
of the boundary element $\sigma_{\ell}$.  The local mesh size is then given as 
$h_{\ell}: = |t_{\ell_{2}} - t_{\ell_{1}}|$ while
$h := \max_{\ell=1,\ldots,N} h_\ell$ is the global mesh size.

\noindent 
For the approximation of the unknown Cauchy datum 
$w = \gamma_{1}^{\textrm{int}} u \in \dualtraceH$ we consider the space 
$S_{h}^{0}(\Sigma) := 
\textrm{span}\left\{\varphi_{\ell}^{0}\right\}_{\ell=1}^{N}$
of piecewise constant basis functions $\varphi_\ell^0$, which is defined
with respect 
to the decomposition $\Sigma_{N}$. The Galerkin-Bubnov variational 
formulation of \eqref{eq_var_form_first_bie} is to find 
$w_{h} \in S_{h}^{0}(\Sigma)$ such that
\begin{equation}
\label{eq_galerkin_var_form_first_bie}
\langle V w_{h}, \tau_{h} \rangle_{\Sigma} = \langle (\frac{1}{2}I + K)g,
\tau_{h} \rangle_{\Sigma} - \langle M_{0}u_{0}, \tau_{h} \rangle_{\Sigma} \quad
\textrm{for all } \tau_{h} \in S_{h}^{0}(\Sigma) \, .
\end{equation}
This is equivalent to the system of linear equations 
$V_{h} \vec{w} = \vec{f}$ where
\begin{equation*}
V_{h}[\ell, k] = \langle V \varphi_{k}^{0}, \varphi_{\ell}^{0} \rangle_{\Sigma},
\quad \vec{f}[\ell] = 
\langle (\frac{1}{2}I + K)g, \varphi_{\ell}^{0}\rangle_{\Sigma}
- \langle M_{0} u_{0}, \varphi_{\ell}^{0} \rangle_{\Sigma},
\quad k,\ell=1,\ldots,N.
\end{equation*}
Due to the ellipticity of the single layer operator $V$ the matrix 
$V_{h}$ is positive definite and therefore the variational formulation
\eqref{eq_galerkin_var_form_first_bie} is uniquely solvable as well. 
Moreover, when assuming $w \in H^{s}(\Sigma)$ for some $s \in [0,1]$, 
there holds the error estimate
\begin{equation*}
\norm{w - w_{h}}_{\dualtraceH} \leq c h^{1/4+s}|w|_{H^{s}(\Sigma)} \, .
\end{equation*}
Using standard arguments we also conclude the error estimate
\begin{equation*}
\norm{w - w_{h}}_{L^{2}(\Sigma)} \leq c h^{s} |w|_{H^{s}(\Sigma)}
\end{equation*}
which implies linear convergence of the $L^{2}(\Sigma)$--error of the 
Galerkin approximation $w_{h}$ if $w \in H^{1}(\Sigma)$ is satisfied. 

\section{Preconditioning strategies}
Since the boundary element discretization is done with respect to the
whole space--time boundary $\Sigma$ we need to have an efficient
iterative solution technique. In fact,
the linear system $V_{h} \vec{w} = \vec{f}$ with the
positive definite but nonsymmetric matrix $V_h$ can be solved by using a
preconditioned GMRES method. Here
we will apply a preconditioning technique 
based on boundary integral operators of opposite order 
\cite{steinbach_wendland}, also known as operator or 
Calderon preconditioning \cite{hiptmair}. Since the single layer
integral operator $V: \dualtraceH \to \traceH$ and
the hypersingular integral operator
$D: \traceH \to \dualtraceH$ are both elliptic, the operator
$DV : \dualtraceH \to \dualtraceH$ behaves like the identity.
Hence we can use the
Galerkin discretization of $D$ as a preconditioner for $V_h$.
But for the Galerkin discretization $D_h$ of the hypersingular integral
operator $D : \traceH \to \dualtraceH$ we need to use a
conforming ansatz space
$Y_{h} = \textrm{span} \left\{\psi_i\right\}_{i=1}^{N} \subset \traceH$
while the discretization of the single layer integral operator $V$
is done with respect to $S_h^0(\Sigma)$. Since the boundary element space $S_h^0(\Sigma)$ of piecewise constant basis functions $\varphi_k^0$ also satisfies $S_h^0(\Sigma) \subset \traceH$ we can choose $Y_h = S_{h}^{0}(\Sigma)$. The inverse hypersingular
operator $D^{-1}$ is spectrally equivalent to the single layer operator
$V$, therefore the approximation of the preconditioning operator corresponds 
to a mixed approximation scheme, and hence we need to assume a 
discrete stability condition to be satisfied.

\begin{theorem}[\cite{hiptmair,steinbach_wendland}]
Assume the discrete stability condition
\begin{equation}
\label{eq_stability_trial_spaces}
\sup_{0 \neq v_{h} \in Y_{h}} \frac{\langle \tau_{h}, v_{h}
  \rangle_{L^{2}(\Sigma)}}{\norm{v_{h}}_{\traceH}} \geq c_{1}^{M}
\norm{\tau_{h}}_{\dualtraceH} \quad \textrm{for all } \tau_{h} \in S_h^0(\Sigma) .
\end{equation}
Then there exists a constant $c_{\kappa} > 1$ such that
\begin{equation*}
\kappa\left(M_{h}^{-1} D_{h} M_{h}^{-\top} V_{h}\right) \leq c_{\kappa}
\end{equation*}
where, for $k,\ell=1,\ldots,N$,
\begin{equation*}
V_{h}[\ell, k] = \langle V \varphi_{k}^0, \varphi_{\ell}^0 \rangle_{\Sigma} \;, 
\quad 
D_{h}[\ell, k] = \langle D \psi_{k}, \psi_{\ell}\rangle_{\Sigma} \;, 
\quad M_{h}[\ell, k] = \langle \varphi_{k}^0, \psi_{\ell} \rangle_{L^{2}(\Sigma)}
\; .
\end{equation*}
\end{theorem}
Thus we can use $C_{V}^{-1} = M_{h}^{-1} D_{h} M_{h}^{-\top}$ as a
preconditioner for $V_h$. Since $M_{h}$ is sparse and spectrally
equivalent to a diagonal matrix, the inverse $M_{h}^{-1}$ can be computed 
efficiently. It remains to define, for given $S_h^0(\Sigma)$,
a suitable boundary element space $Y_h$ such that the
stability condition \eqref{eq_stability_trial_spaces} is satisfied.
In what follows we will discuss a possible choice.

If we choose $Y_{h} = S_{h}^{0}(\Sigma)$ for the discretization of 
the hypersingular operator $D$, then $M_h$ becomes diagonal
and is therefore easily invertible.
In order to prove the stability condition 
\eqref{eq_stability_trial_spaces} we need to establish 
the $\traceH$--stability of the $L^{2}(\Sigma)$--projection 
$Q^{0}_{h} : L^{2}(\Sigma) \to S_{h}^{0}(\Sigma) \subset L^2(\Sigma)$ 
which is defined as
\begin{equation*}
\langle Q_{h}^{0}v, \tau_{h} \rangle_{L^{2}(\Sigma)} = 
\langle v, \tau_{h} \rangle_{L^{2}(\Sigma)} \quad \textrm{for all } 
\tau_{h} \in S_{h}^{0}(\Sigma).
\end{equation*}
Following \cite{steinbach_stability_l2}, and when assuming 
local quasi-uniformity of the boundary element mesh $\Sigma_N$ we
are able to establish the stability of $Q_{h}^{0}:\traceH \rightarrow \traceH$,
see \cite{dohr_niino_steinbach} for a more detailed discussion:
For $\ell=1,...,N$ we define $I(\ell)$ to be the index set of the boundary 
element $\sigma_{\ell}$ and all its adjacent elements. We assume the boundary 
element mesh $\Sigma_{N}$ to be locally quasi-uniform, i.e. there 
exists a constant $c_{L} \geq 1$ such that
\begin{equation*}
\frac{1}{c_{L}} \leq \frac{h_{\ell}}{h_{k}} \leq c_{L} \quad 
\textrm{for all } k \in I(\ell) \textrm{ and  } \ell=1,...,N.
\end{equation*}
In this case the operator $Q_{h}^{0}:\traceH \rightarrow \traceH$ is 
bounded, i.e. there exists a constant $c_{S}^{0} > 0$ such that
\begin{equation}
\label{eq_stability_l2_h14_piecewise_constant}
\norm{Q_{h}^{0}v}_{\traceH} \leq c_{S}^{0} \norm{v}_{\traceH} \quad 
\textrm{for all } v \in \traceH.
\end{equation}
By using the stability estimate 
\eqref{eq_stability_l2_h14_piecewise_constant} we can conclude
\begin{equation*}
\frac{1}{c_{S}^{0}} \norm{\tau_{h}}_{\dualtraceH} \leq \sup_{0 \neq v_{h} \in S_{h}^{0}(\Sigma)} \frac{\langle \tau_{h}, v_{h}\rangle_{L^{2}(\Sigma)} }{\norm{v_{h}}_{\traceH}} \quad \textrm{for all } \tau_{h} \in S_{h}^{0}(\Sigma).
\end{equation*}
Hence the stability condition \eqref{eq_stability_trial_spaces} holds and 
we can use $C_{V}^{-1} = M_{h}^{-1} D_{h} M_{h}^{-\top}$ as a
preconditioner for $V_h$.

\section{Numerical results}
For the numerical experiments we choose $\Omega =(0,1)$, $T=1$, and we 
consider the model problem \eqref{eq_model_problem} with homogeneous 
Dirichlet conditions $g=0$, and some given initial datum 
$u_{0}$ satisfying the compatibility conditions $u_{0}(0) = u_{0}(1) = 0$. 
The Galerkin boundary element discretization of the variational formulation
\eqref{eq_var_form_first_bie} is done by piecewise constant basis functions.
The resulting system of linear equations $V_{h}\vec{w} = \vec{f}$ is solved 
by using the GMRES method. As a preconditioner we use the discretization 
$C_{V}^{-1} = M_{h}^{-1} D_{h} M_{h}^{-\top}$
of the hypersingular operator $D$ with piecewise constant basis functions. 

\paragraph{Uniform refinement}
The first example corresponds to the initial datum 
$u_{0}(x) = \sin \, 2 \pi x$ and a globally uniform boundary element mesh
of mesh size $h=2^{-L}$.
Table \ref{tab_prec_uniform} shows the $L^{2}(\Sigma)$--error 
$\| w-w_h \|_{L_2(\Sigma)}$ and the 
estimated order of convergence (eoc), which is linear as expected.
Moreover, the condition numbers of the stiffness matrix $V_h$ and of the
preconditioned matrix $C_V^{-1}V_h$ 
as well as the number of iterations to reach a relative accuracy of
{\bf $10^{-8}$} are given which confirm the theoretical estimates.

\begin{table}[htb]
\caption{Error, condition and iteration numbers in the case of 
uniform refinement}
\label{tab_prec_uniform} 
\begin{tabular}{p{0.5cm}p{1cm}p{2cm}p{1cm}p{1cm}p{0.5cm}p{1.5cm}p{0.5cm}}
\hline\noalign{\smallskip}
L & $N$ & $\norm{w - w_{h}}_{L_{2}(\Sigma)}$ & eoc & $\kappa(V_{h})$ & It. & $\kappa(C_{V}^{-1} V_{h})$ & It.  \\
\noalign{\smallskip}\svhline\noalign{\smallskip}
0&2&2.249&-&1.001&1&1.002&1\\
1&4&1.311&0.778&2.808&2&1.279&2\\
2&8&0.658&0.996&4.905&4&1.422&4\\
3&16&0.324&1.021&7.548&8&1.486&8\\
4&32&0.160&1.017&11.140&16&1.541&14\\
5&64&0.079&1.010&16.724&31&1.563&13\\
6&128&0.040&1.006&13.470&41&1.590&13\\
7&256&0.020&1.003&22.053&50&1.615&12\\
8&512&0.010&1.001&32.043&59&1.636&12\\
9&1024&0.005&1.001&60.957&70&1.777&11\\
10&2048&0.002&1.000&88.488&82&1.762&11\\
11&4096&0.001&1.000&125.957&96&1.765&10\\
\noalign{\smallskip}\hline\noalign{\smallskip}
\end{tabular}
\end{table}

\paragraph{Adaptive refinement}
For the second example we consider the initial datum 
$u_{0}(x) = 5 e^{-10 x} \sin \, \pi x$ which motivates the use of 
a locally quasi-uniform boundary element mesh resulting from some
adaptive refinement strategy. The numerical results as given in
Table \ref{tab_prec_adaptive} again confirm the theoretical findings,
in particular the robustness of the proposed preconditioning strategy
in the case of an adaptive refinement which is not the case
when using none or only diagonal preconditioning
$\widetilde{C}_{V} = \mathrm{diag}V_{h}$.

\begin{table}[htb]
\caption{Error, condition and iteration numbers in the case of 
adaptive refinement}
\label{tab_prec_adaptive}
\begin{tabular}{p{0.5cm}p{1cm}p{2cm}p{1.3cm}p{0.7cm}p{1.2cm}p{0.5cm}p{1.2cm}p{0.5cm}}
\hline\noalign{\smallskip}
L & $N$ & $\norm{w - w_{h}}_{L_{2}(\Sigma)}$ & $\kappa(V_{h})$ & It. & 
$\kappa(\widetilde{C}_V^{-1} V_{h})$& It. & $\kappa(C_{V}^{-1} V_{h})$ & It.\\
\noalign{\smallskip}\svhline\noalign{\smallskip}
0&2&1.886&1.00&2&1.001&2&1.002&2\\
1&3&1.637&3.97&3&2.553&3&1.16&3\\
2&5&1.272&12.23&5&4.055&4&1.166&4\\
3&7&0.914&34.21&7&3.611&6&1.156&6\\
4&9&0.615&92.08&9&3.164&8&1.149&8\\
5&11&0.401&118.59&11&2.945&10&1.224&10\\
6&13&0.267&338.26&13&2.803&12&1.21&12\\
7&20&0.166&621.77&20&3.524&18&1.197&13\\
8&31&0.101&1608.08&31&4.457&27&1.252&12\\
9&47&0.063&2344.90&47&5.779&32&1.574&11\\
10&74&0.039&6141.47&74&8.348&37&1.692&11\\
11&114&0.024&8409.92&114&10.950&42&1.561&10\\
12&177&0.015&23007.60&173&14.324&47&1.716&10\\
13&278&0.010&27528.30&200&21.094&53&1.677&10\\
\noalign{\smallskip}\hline\noalign{\smallskip}
\end{tabular}
\end{table}

\section{Conclusions and outlook}
In this note we have described a space--time boundary element
discretization of the spatially one-dimensional heat equation
and an efficient and robust preconditioning strategy which is based on
the use of boundary integral operators of opposite orders, but which
requires a suitable stability condition for the boundary element
spaces used for the discretization.
In the particular case of the spatially one-dimensional heat equation
we can use the space $S_h^0(\Sigma)$ of piecewise constant basis functions 
to discretize both the single layer and the hypersingular boundary
integral operator $V$ and $D$, respectively. This is due to the
inclusion $S_h^0(\Sigma) \subset \traceH$ where the latter is
the Dirichlet trace space of the anisotropic Sobolev space
$H^{1,1/2}(Q)$. In the case of a spatially two- or three-dimensional domain
$\Omega$ a conformal approximation of the Dirichlet trace space
$H^{1/2,1/4}(\Sigma)$ and therefore the discretization of the hypersingular
integral operator $D$ requires the use of continuous basis functions.
Hence, to ensure the stability condition
\eqref{eq_stability_trial_spaces} we may use the space $S_h^1(\Sigma)$
of piecewise linear and continuous basis functions for the discretization
of $V$ and $D$, respectively, 
see \cite[Theorem 3.2]{steinbach_stability_l2}, and when assuming
some appropriate mesh conditions locally  
\cite[Section 4]{steinbach_stability_l2}. However, due to the approximation
properties of $S_h^1(\Sigma)$ such an approach is restricted to spatial
domains $\Omega$ with smooth boundary where the unknown flux is continuous.

\begin{figure}[ht]
\sidecaption[t]
\includegraphics[]{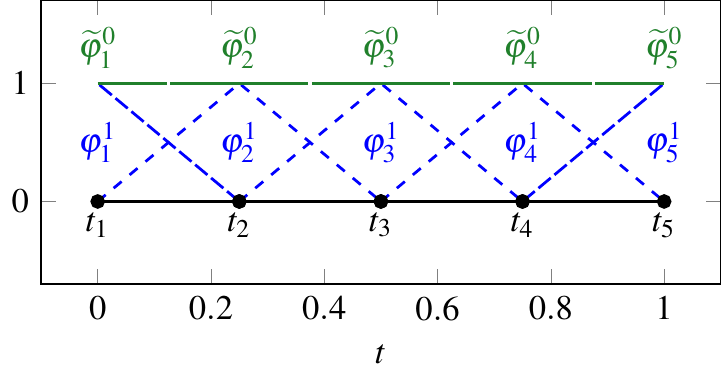}
\caption{Sample dual mesh. The piecewise linear and continuous functions $\varphi_{i}^{1}$ are used for the discretization of $D$. The piecewise constant basis functions $\tilde{\varphi}_{i}^{0}$ are used for the discretization of $V$}
\label{fig_dual_mesh_1d}
\end{figure}

When using the discontinuous boundary element space $S_h^0(\Sigma)$ for the
approximation of the unknown flux we need to choose an appropriate
boundary element space $Y_h$ to ensure the stability condition
\eqref{eq_stability_trial_spaces}. A possible approach is the use of a
dual mesh using piecewise constant basis functions 
for the approximation of $V$, and piecewise linear and continuous basis 
functions for the approximation of $D$, see
Fig. \ref{fig_dual_mesh_1d} for the situation in 1D.
For a more detailed analysis of the proposed preconditioning
strategy and suitable choices of stable boundary element spaces
we refer to \cite{dohr_niino_steinbach}.

\noindent
An efficient solution of local Dirichlet boundary value problems is an
important tool when considering domain decomposition methods for the
heat equation, see e.g. \cite{steinbach_wendland_97} in the case of 
the Laplace equation. Moreover, the preconditioning strategy of
using operators of opposite order can also be used when considering
related Schur complement systems on the skeleton, as they also
appear in tearing and interconnecting domain decomposition methods, 
see, e.g., \cite{langer_steinbach}. This also covers the coupling
of space--time finite and boundary element methods. Related results
on the stability and error analysis as well as on efficient
solution strategies for space--time domain decomposition methods
will be published elsewhere.

\acknowledgement{This work was supported by the International Research Training Group 1754, funded by the German Research Foundation (DFG) and the Austrian Science Fund (FWF). Additionally, S.~Dohr would like to acknowledge the financial support provided by the University of Bergen.}

\input{referenc}

\end{document}

%% file: referenc.tex
%
%
%

%% file: main.bbl
\begin{thebibliography}{99.}%
%
%
\bibitem{costabel} Costabel, M.: Boundary integral operators for the heat 
equation. Integral Equations Operator Theory, 13: 498--552, 1990.

\bibitem{dohr_niino_steinbach}
Dohr, S., Niino, K., Steinbach, O.:
Preconditioned space--time boundary element methods for the heat
equation, in preparation.

\bibitem{hiptmair} Hiptmair, R.: Operator Preconditioning. 
Comput. Math. Appl., 52: 699--706, 2006.

\bibitem{langer_steinbach} Langer, U., Steinbach, O.: Boundary element tearing
  and interconnecting methods. Computing, 71: 205--228, 2003. 

\bibitem{lions_magenes_2} Lions, J.L., Magenes, E.: Non-Homogeneous Boundary 
Value Problems and Applications II. Springer, Berlin-Heidelberg, 1972.
  
\bibitem{noon} Noon, P.J.: The Single Layer Heat Potential and Galerkin 
Boundary Element Methods for the Heat Equation. PhD thesis, University
of Maryland, 1988.

\bibitem{steinbach_stability_l2}
Steinbach, O.: On the stability of the $L_2$ projection in fractional
Sobolev spaces. Numer. Math., 88: 367--379, 2001.

\bibitem{steinbach}
Steinbach, O.: Space--time finite element methods for parabolic problems. 
Comput. Meth. Appl. Math., 15: 551-566, 2015. 

\bibitem{steinbach_wendland_97} Steinbach, O., Wendland, W.L.: Efficient 
preconditioners for boundary element methods and their use in 
domain decomposition. In: Domain Decomposition Methods in Sciences in
Engineering: 8th International Conference, Beijing, China 
(R. Glowinski et. al. eds.), John Wiley, pp. 3--18, 1997. 

\bibitem{steinbach_wendland} Steinbach, O., Wendland, W.L.: 
The construction of some efficient preconditioners in the boundary element 
method. Adv. Comput. Math., 9: 191--216, 1998.

\end{thebibliography}
